\documentclass[12pt]{amsart}
\usepackage{graphicx}
\usepackage{amsmath}
\usepackage{amsfonts}
\usepackage{amssymb}

\newtheorem{theorem}{Theorem}

\newtheorem{corollary}[theorem]{Corollary}

\newtheorem{proposition}[theorem]{Proposition}
\newtheorem{remark}{Remark}

\def\F{\Phi}
\def\H{\mathcal H}
\def\aa{\alpha}

\def\g{\gamma}
\def\d{\delta}

\def\SS{{\bf S}}
\def\Z{{\mathbb Z}}
\def\R{{\mathbb R}}

\def\o{\omega}

\def\ps{\varphi}
\def\a{\alpha}
\def\b{\beta}
\def\bb{\gamma}

\def\wa{\widehat\alpha}
\def\la{\bar\alpha}
\def\lb{\bar\beta}

\def\f{\varphi}
\def\p{\pi}
\def\ee{\varepsilon}

\def\l{\lambda}

\begin{document}

\title[Strongly-cyclic branched coverings of knots]{Strongly-cyclic branched coverings of knots via $(g,1)$-decompositions}

\thanks{Work performed under the auspices of G.N.S.A.G.A.
of C.N.R. of Italy and supported by M.U.R.S.T., by the University
of Bologna, funds for selected research topics, and by the Russian
Foundation for the Basic Researches. The third named author was
supported by the INTAS project ``CalcoMet-GT' 03-51-3663 and by
RFRB}


\author[Cristofori]{Paola Cristofori}
\address{Department of Mathematics, University of Modena and Reggio Emilia, Italy}
\email{cristofori.paola@unimo.it}
\author[Mulazzani]{Michele Mulazzani}
\address{Department of Mathematics and C.I.R.A.M., University of Bologna, Italy}
\email{mulazza@dm.unibo.it}
\author[Vesnin]{Andrei Vesnin}
\address{Sobolev Institute of Mathematics, Novosibirsk 630090, Russia}
\email{vesnin@math.nsc.ru}

\subjclass{Primary 57M12, 57M25; Secondary 20F05, 57M05, 20F10}
\keywords{$(g,1)$-knots, $(g,1)$-decompositions, cyclic branched
coverings, presentations of groups, Heegaard splittings, Heegaard
diagrams.}

\date{}

\maketitle

\begin{abstract}
Strongly-cyclic branched coverings of knots are studied by using
their $(g,1)$-decompositions. Necessary and sufficient conditions
for the existence and uniqueness of such coverings are obtained.
It is also shown that their fundamental groups admit geometric
$g$-words cyclic presentations.
\end{abstract}

\section{Introduction}

A $(g,1)$-decomposition of a knot in a closed orientable
3-manifold is a genus $g$ Heegaard splitting of the manifold, such
that the knot is trivially split in two arcs by the Heegard
surface (a more detailed definition will be given in Section~2).
Any knot $K$ in a 3-manifold $N$ admits a $(g,1)$-decomposition,
for a certain $g\geq h(N)$, where $h(N)$ is the Heegaard genus of
$N$. In the following a knot admitting a $(g,1)$-decomposition
will be called a $(g,1)$-knot. The particular case of
$(1,1)$-knots has recently been investigated in several papers
from different points of view (see references in \cite{CM1}).

In this paper we study strongly-cyclic branched coverings of
knots, by using topological and algebraic properties of their
$(g,1)$-decompositions. Some of the obtained results naturally
gene\-ra\-lize results from \cite{CM0} and \cite{M} concerning
$(1,1)$-knots.

In Section~2, we prove a representation theorem of $(g,1)$-knots
in terms of Heegaard diagrams and obtain a presentation for the
fundamental group of the knot complement. In Section~3, we
consider strongly-cyclic branched coverings of knots and give
conditions for their existence and uniqueness, using
$(g,1)$-decompositions. As a relevant example we consider, in
Section~4, the case of generalized periodic Takahashi manifolds.
Section~5 is devoted to fundamental groups of strongly-cyclic
branched coverings of knots. These groups admit geometric $g$-word
cyclic presentations, a straightforward generalization of cyclic
presentations of groups. Furthermore, we describe an algorithm to
obtain such a presentation, starting from the group of the knot
complement.

For the theory of Heegaard splittings of 3-manifolds we refer to
\cite{He}.

\section{$(g,1)$-decompositions of knots}

The notion of $(g,b)$-decomposition of a link $L$ in a closed
orientable 3-manifold has been introduced in \cite{Do}, as a
generalization of the classical bridge (or plat) decomposition of
links in $\SS^3$.

A non-empty set $\{a_1,\ldots ,a_b\}$ of mutually disjoint arcs
properly embedded in a genus $g$ handlebody $\H_g$ is called
\textit{trivial} if there exist $b$ mutually disjoint discs
$D_1,\ldots ,D_b\subset \H_g$ such that $a_i\cap D_i=
a_i\cap\partial D_i=a_i$, $a_i\cap D_j=\emptyset$ and $\partial
D_i-a_i\subset\partial \H_g$ for all $i,j=1,\ldots ,b$ and $i\neq
j$.

Let $\H_g$ and $\H_g'$ be the two handlebodies of a genus $g$
Heegaard splitting of a closed orientable 3-manifold $N$ and let
$S_g=\partial \H_g=\partial \H_g'$. A link $L\subset N$ is said to
be in {\it $b$-bridge position\/} with respect to $S_g$ if: (i)
$L$ intersects $S_g$ transversally and (ii) $L\cap \H_g$ and
$L\cap \H_g'$ are both sets of $b$ mutually disjoint properly
embedded trivial arcs. This splitting is called a {\it
$(g,b)$-decomposition of $L$\/}.

Equivalently, a $(g,b)$-decomposition of $L$ can be obtained by
taking two genus $g$ handlebodies $\H_g,\H_g'$, two trivial sets
of arcs \hbox{$A_b\subset \H_g,A_b'\subset \H_g'$} and an
attaching homeomorphism $\f:(\partial \H_g', \partial
A_b')\to(\partial \H_g,\partial A_b)$ such that
$$(N,L)\cong (\H_g,A_b)\cup_{\f}(\H_g',A_b').$$

Note that a $(0,b)$-decomposition is the classical bridge
decomposition for links in a 3-sphere.

A link $L$ is called a {\it $(g,b)$-link\/} if it admits a
$(g,b)$-decomposition. Of course, a $(g,1)$-link is a knot, for
every $g\geq 0$, and a $(0,1)$-knot is a trivial knot in $\SS^3$.
Obviously, a $(g,b)$-link is also a $(g',b')$-link for any $g'\ge
g$ and $b'\ge b$.

\begin{figure}[h]
\begin{center}
\includegraphics*[totalheight=2cm]{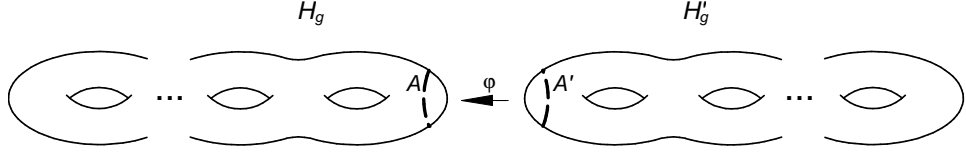}
\end{center}
\caption{A $(g,1)$-decomposition} \label{decomp}
\end{figure}

In the following we shall restrict our attention to $(g,1)$-knots.
A $(g,1)$-decomposition of a knot is represented in Figure
\ref{decomp}.

A knot $K\subset N$ is called a {\it $(g,0)$-knot\/} if it can be
embedded in a genus $g$ Heegaard surface of $N$. It is easy to see
that a $(g,0)$-knot is also a $(g,1)$-knot (just push an arc of
$K$ into a handlebody of the splitting and the remaining part of
the knot into the other handlebody). Moreover, when $b>0$, a
$(g,b)$-knot is also a $(g+1,b-1)$-knot.

\begin{remark}\label{Remark 1}
{\rm Let $K$ be a knot in a 3-manifold $N$ of Heegaard genus $h$.
It is easy to see that, for any integer $b\geq 0$, there exists an
integer $g\geq h$ such that $K$ is a $(g,b)$-knot in $N$.
Analogously, for any integer $g\geq h$ there exists an integer
$b\geq 0$ such that $K$ is a $(g,b)$-knot in $N$.}
\end{remark}

The following presentation theorem shows that a $(g,1)$-knot can
be described by the images of a certain set of meridian curves.

\begin{theorem} \label{basic}
Let $K\subset N$ be a knot with a $(g,1)$-decomposition
$(N,K)=(\H_g,A)\cup_{\f}(\H_g',A').$ Then:
\begin{itemize}
\item [(i)] the knot $K$ is completely determined, up to
equivalence\footnote{Two knots $K\subset N$ and $K'\subset N'$ are
considered {\it equivalent} if there exists a homeomorphism
$\phi:N\to N'$ such that $\phi(K)=K'$.}, by
$\f(\beta_1'),\f(\beta_2'),\ldots,\f(\beta_g')$, where
$\beta_1'=\partial E_1'$, $\ldots$ , $\beta_g'=\partial E_g'$, and
$E_1',\ldots,E_g'$ is a complete system of meridian disks of
$\H_g'$ not intersecting $A'$; \item [(ii)] if $(N,\bar
K)=(\H_g,A)\cup_{\bar{\f}}(\H_g',A')$ is a $(g,1)$-decomposition
of a knot $\bar K$ such that $\bar{\f}(\beta_j')$ is isotopic to
$\f(\beta_j')$ in $\partial \H_g-\partial A$, for $j=1,\ldots,g$,
then $\bar K$ is equivalent to $K$; \item [(iii)] if
$A''\subset\partial \H_g$ is an arc such that $\partial
A''=\partial A$ and $A''\cap\f (\beta_j')=\emptyset$ for
$j=1,\ldots,g$, then $A\cup A''$ is a $(g,1)$-knot equivalent to
$K.$
\end{itemize}
\end{theorem}
\begin{proof}
(i) It follows from the fact that two properly embedded trivial
arcs in a ball $B$, with the same endpoints,  are isotopic rel
$\partial B$. Statements (ii) and (iii) are straightforward.
\end{proof}


Now we find a presentation for the fundamental group of
$(g,1)$-knots complements.

Let $K$ be a knot in a manifold $N$ with a $(g,1)$-decomposition
$(N,K)=(\H_g,A)\cup_{\f}(\H_g',A')$. Moreover, let $\tau
:(\H_g,A)\to (\H_g',A')$ be a fixed homeomorphism.

\begin{figure}[h]
\begin{center}
\includegraphics*[totalheight=2.8cm]{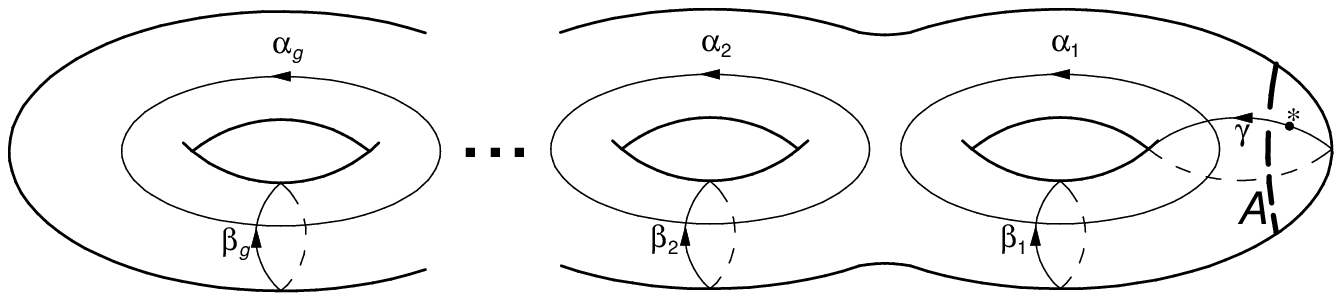}
\end{center}
\caption{} \label{generat}
\end{figure}

Consider the loops $\a_1,\ldots,\a_g,\b_1,\ldots,\b_g,\bb$ in
$\partial \H_g$, as depicted in Figure \ref{generat}. Fix a base
point $*$ on the curve $\gamma$ and, for each $i=1,\ldots,g$,
define the loops $\la _i = \xi_i \cdot \a _i \cdot \xi_i^{-1}$ and
$\lb _i=\eta_i\cdot\b _i\cdot\eta_i^{-1}$, where $\xi _i$ and
$\eta _i$ are paths connecting $*$ to $\a _i$ and $\b _i$,
respectively. Obviously, the cycles $\la_i$ and $\lb_i$ are
homologous to $\a_i$ and $\b_i$, respectively, for each
$i=1,\ldots,g$. Moreover, we set $\la'_i=\tau(\la_i)$,
$\lb'_i=\tau(\lb_i)$, for $i=1,\ldots,g$, $\gamma'=\tau(\gamma)$
and $*'=\tau(*)$. Then $*'=\ps^{-1}(*)$, up to isotopy. The
homotopy classes of $\la _1,\ldots,\la _g,\lb _1,\ldots,\lb _g$
and $\gamma$ generate $\p_1(\partial \H_g-\partial A,*)$ and the
homotopy classes of $\la' _1,\ldots,\la' _g,\lb' _1,\ldots,\lb'
_g$ and $\gamma'$ generate $\p_1(\partial \H_g'-\partial A',*')$.
In order to simplify the notation, in the following we use the
same symbol for loops (resp. for cycles) and for corresponding
homotopy classes (resp. homology classes).


\begin{proposition} \label{fundamental} The fundamental group of
the complement of a knot $K\subset N$, with $(g,1)$-decomposition
$(N,K)=(\H_g,A)\cup_{\f}(\H_g',A')$ admits the presentation
\begin{equation} \label{presentation}
\p_1(N-K,*) = \langle
\la_1,\ldots,\la_g,\gamma\,|\,r_i(\la_1,\ldots,\la_g,\gamma),
\quad i=1,\ldots,g \rangle,
\end{equation}
where $r_i(\la_1,\ldots,\la_g,\gamma)=\iota_\#\f_\#(\lb_i')$,
being $\f_\#$ and $\iota_\#$ the homomorphisms of fundamental
groups induced by the attaching homeomorphism $\f : (\partial
\H_g',
\partial A') \to (\partial \H_g, \partial A)$ and the inclusion \hbox{$\iota : (\partial \H_g,\partial A)\to
(\H_g,A)$,}
respectively.
\end{proposition}
\begin{proof}
 We have $\pi_1(\H_g-A,*)= \langle \la_1,\ldots,\la_g,\lb_1,\ldots,\lb_g,\gamma\, |
 \,\lb_1,\ldots,\lb_g \rangle$ and
$\pi_1(\H_g'-A',*') = \langle
\la'_1,\ldots,\la'_g,\lb'_1,\ldots,\lb'_g,\gamma'\, |
\,\lb'_1,\ldots,\lb'_g \rangle$. Applying the Seifert-van~Kampen
theorem to the $(g,1)$-decomposition of $K$, we get
 \begin{eqnarray*}
\p_1(N-K,*) =  \langle \la_1, \ldots, \la_g, \lb_1, \ldots,
\lb_g, \gamma, \la'_1, \ldots, \la '_g, \lb'_1, \ldots, \lb'_g, \gamma'\, | \\
 \lb_i,\lb'_i,\ps(\la'_i)\la_i^{'-1}, \ps(\lb'_i)\lb_i^{'-1},
\ps(\gamma')\gamma^{'-1}, \quad  i=1,\ldots,g \rangle \\
 = \langle \la_1, \ldots, \la_g, \lb_1, \ldots, \lb_g, \gamma \, |
\, \lb_i, \f(\lb_i'), \quad  i=1,\ldots,g \rangle \\
= \langle \la_1, \ldots, \la_g, \gamma \, | \,
r_i(\la_1,\ldots,\la_g,\gamma), \quad  i=1,\ldots,g \rangle,
 \end{eqnarray*}
where $r_i(\la_1,\ldots,\la_g,\gamma)$ is obtained by erasing all
$\lb_j$'s terms from $\f (\lb_i')$. Since $\pi_1(\partial
\H_g-\partial A,*) = \langle \la_1, \ldots, \la_g, \lb_1, \ldots,
\lb_g, \gamma \, | \,\emptyset \rangle$, the statement is proved.
\end{proof}

As a consequence of the above proposition we have the following:

\begin{corollary}\label{homology} The first homology group of $N-K$ admits the
following presentation:
\begin{equation} \label{abelian}
H_1(N-K) =  \langle \a_1, \ldots, \a_g, \bb \, | \, a_{i1} \a_1 +
\ldots + a_{ig} \a_g + b_i \bb, \ i=1, \ldots, g \rangle
\end{equation}
where, for each $i,j=1,\ldots,g$, $a_{ij}$ (resp. $b_i$) is the
exponent sum of $\la_j$ (resp. $\gamma$) in the word $r_i$ of the
presentation {\rm (\ref{presentation})}.
\end{corollary}

\section{Strongly-cyclic branched coverings of knots}

Let $K\subset N$ be a knot in a 3-manifold $N$ and $n>1$. Then an
$n$-fold cyclic covering $f : M\to N$ branched over $K$ is called
{\it strongly-cyclic\/} if the branching index of $K$ is $n$. This
means that the fiber $f^{-1}(x)$ of each point $x\in K$ contains a
single point. In this case the homology class of a meridian loop
$\mu$ around $K$ is mapped by the associated monodromy
\hbox{$\o_{f}:H_1(N-K)\to\Z_n$} to a generator of the cyclic group
$\Z_n$, and up to equivalence  we can always suppose
$\o_{f}(\mu)=1$. We recall that two $n$-fold cyclic coverings  $f'
: M'\to N$ and $f'' : M''\to N$, branched over $K\subset N$, with
monodromies $\o_{f'}$ and $\o_{f''}$ respectively, are equivalent
if and only if there exists $u\in\Z_n$, with $\gcd(u,n)=1$, such
that $\o_{f''}=u\o_{f'}$, where $u\o_{f'}$ is the multiplication
of $\o_{f'}$ by $u$. For knots in $\SS^3$, strongly-cyclic
branched coverings and cyclic branched coverings are equivalent
notions. The same property occurs in the general case when $n$ is
prime.

Now we discuss strongly-cyclic branched coverings of knots via
their $(g,1)$-decompositions. In particular, we give existence
conditions and compute their number, up to equivalence.

Let $K$ be a knot in a 3-manifold $N$, and let
\begin{equation} \label{homology}
H_1(N)=\Z^d\oplus\Z_{t_1}\oplus\Z_{t_2}\oplus\cdots\oplus\Z_{t_s},
\end{equation}
with $0\le d$ and $1<t_1\,\vert\,t_2\,\vert\,\cdots\,\vert\, t_s$,
be the canonical decomposition of the first integer homology group
of $N$.

Let us fix a $(g,1)$-decomposition for $K$. Referring to
presentation (2), let $H\in {\mathcal M}_g(\Z)$ be the matrix with
coefficients $h_{ij}=a_{ij}$ and $H'\in {\mathcal M}_{g,g+1}(\Z)$
the augmented matrix with coefficients $h_{ij}'=a_{ij}$, for
$i,j=1,\ldots,g$ and $h_{ij}'=b_i$ for $j=g+1$ and $i=1,\ldots,g$.
Let $e_1,\ldots,e_r$ (resp. $e'_1,\ldots,e'_{r'}$) be the
invariant factors of $H$ (resp. $H'$). Moreover, set
$e_{r+1}=\ldots=e_g=0$ and $e'_{r'+1}=\ldots=e'_g=0$.

Observe that $H$ is a matrix presentation for $H_1(N)$. Therefore
$s\le r$ and we have: $d=g-r$, $e_i=1$, for $1\le i\le r-s$, and
$e_i=t_{i+s-r}$, for $r-s<i\le r$. Moreover, $H_1(N)$ is finite if
and only if $\det H\ne 0$. In this case $\# H_1(N)=\vert\det
H\vert$.

With the above notations, we have the following result.

\begin{theorem} \label{existence}
Let $K$ be a $(g,1)$-knot in a 3-manifold $N$ with first homology
group given by (\ref{homology}). Then $K$ admits $n$-fold
strongly-cyclic branched coverings if and only if $\gcd
(e_i,n)=\gcd (e'_i,n)$, for each $i=1,\ldots,g$. In this case the
number of non-equivalent $n$-fold strongly-cyclic branched
coverings of $K$ is $C_{N,n}=n^d\gcd (t_1,n)\gcd (t_2,n)\cdots\gcd
(t_s,n)$, which does not depend on $K$.
\end{theorem}
\begin{proof}
Since we suppose $\o_f(\gamma)=1$, an $n$-fold strongly-cyclic
branched covering of $K$ is determined by a solution of the linear
system (in $\Z_n$) $H{\bf x}=-{\bf b}$, and different solutions
define non-equivalent coverings. So the result comes from
\cite{BS}, where systems of linear congruences are investigated.
\end{proof}

As a consequence, we state a uniqueness result.

\begin{corollary} \label{uniqueness}
A knot $K$ in a 3-manifold $N$ admits a unique $n$-fold
strongly-cyclic branched covering, up to equivalence, if and only
if $H_1(N)$ is finite and $\gcd(\#H_1(N),n)=1$.
\end{corollary}
\begin{proof} Let $H$ be the matrix
associated to a $(g,1)$-decomposition of $K$, according to the
previous notations. If $H_1(N)$ is finite and $\gcd(\#
H_1(N),n)=1$, then $H$ admits inverse in $\Z_n$, since in this
case $\vert\det H\vert=\#H_1(N)\ne 0$. Therefore, ${\bf
x}=-H^{-1}{\bf b}$ defines the unique $n$-fold strongly-cyclic
branched covering of $K$. On the other hand, if $K$ admits a
unique $n$-fold strongly-cyclic branched covering, then
$\gcd(e_i,n)=1$ for $i=1,\ldots,g$. Therefore $e_i\ne 0$, for each
$i=1,\ldots,g$. As a consequence, $\vert\det H\vert=e_1e_2\cdots
e_g\ne 0$ and $\gcd(\vert\det H\vert,n)=1$. Since in this case
$\#H_1(N)=\vert\det H\vert$, the statement is proved.
\end{proof}

The previous corollary holds, for any $n$, for knots in homology
spheres (possibly $\SS^3$), since in this case the first homology
group is trivial.

\medskip

\noindent {\bf Examples.} We denote by $M_k$ the connected sum of
$k$ copies of $\SS^2\times\SS^1$, then $H_1(M_k)=\Z^k$.
\begin{itemize} \item The trivial knot $T$ in $M_g$ is a
$(g,1)$-knot defined by $\f(\beta_i')=\beta_i$, for each
$i=1,\ldots,g$. So $H$ and $H'$ are both null matrices, and $T$
admits exactly $n^g$ non-equivalent $n$-fold strongly-cyclic
branched coverings, each of them yielding $M_{ng}$. \item The
``core'' knot $C$ in $M_g$ is defined by $\f(\beta_1')=\gamma$,
and $\f(\beta_i')=\beta_i$, for each $i=2,\ldots,g$. In this case
$H$ is the null matrix but $H'$ has $e'_1=1$, and therefore $C$
admits no $n$-fold strongly-cyclic branched covering, for any
$n>1$.
\end{itemize}

\section{The case of generalized periodic Takahashi manifolds}

\begin{figure}[h]
\begin{center}
\includegraphics*[totalheight=4.5cm]{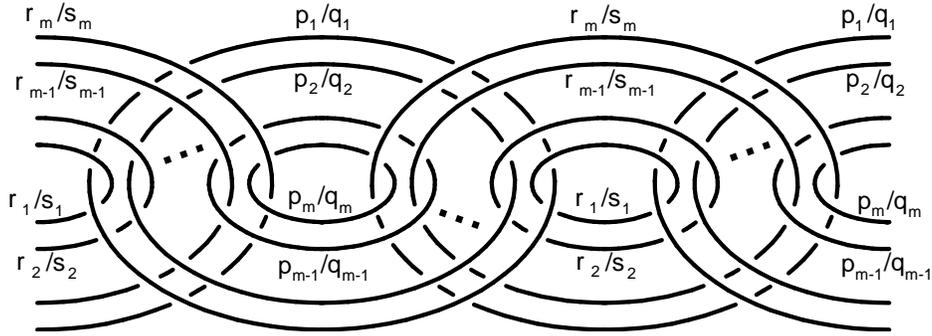}
\end{center}
\caption{$T_{n,m} (p_1/q_1, \ldots, p_m/q_m ; r_1/s_1, \ldots,
r_m/s_m)$} \label{osaka}
\end{figure}

Generalized periodic Takahashi manifolds   were constructed in
\cite{MV} by Dehn surgery on a $n$-periodic $2mn$-component link,
as illustrated in Figure \ref{osaka}, and denoted by $T_{n,m}
(p_1/q_1, \ldots, p_m/q_m ; r_1/s_1, \ldots, r_m/s_m)$, with
$m,n>0$, $p_i,q_i,r_i,s_i\in\mathbb Z$, $\gcd(p_i,q_i) =
\gcd(r_i,s_i)=1$, for each $i=1,\ldots,m$. For $m=1$, we obtain
the periodic Takahashi manifolds $T_{n}(p/q,r/s)$ introduced in
\cite{Ta}.

In \cite{MV} it is proved that generalized periodic Takahashi
manifolds are strongly-cyclic branched coverings of certain knots
in connected sums of lens spaces. That result is improved by the
following.

\begin{proposition}\label{Takahashi}
For each $n>1$, the generalized periodic Takahashi manifold
$T_{n,m} (p_1/q_1,\ldots,p_m/q_m;r_1/s_1,\ldots,r_m/s_m)$ is an
$n$-fold strongly-cyclic branched covering of a $(2m,0)$-knot (and
therefore of a $(2m,1)$-knot).
\end{proposition}
\begin{proof} By
\cite[Theorem 8]{MV},
$T_{n,m}(p_1/q_1,\ldots,p_m/q_m;r_1/s_1,\ldots,r_m/s_m)$ is an
$n$-fold strongly-cyclic covering of the connected sum of lens
spaces $N=L(p_1,q_1)\#L(r_1,s_1)\#\cdots\#L(p_m,q_m)\#L(r_m,s_m)$,
obtained by rational surgery on a $2m$-component trivial link $L$,
branched over a knot $K\subset N$, as shown in Figure~\ref{Fig.
0}.

\begin{figure}[h]
\begin{center}
\includegraphics*[totalheight=6cm]{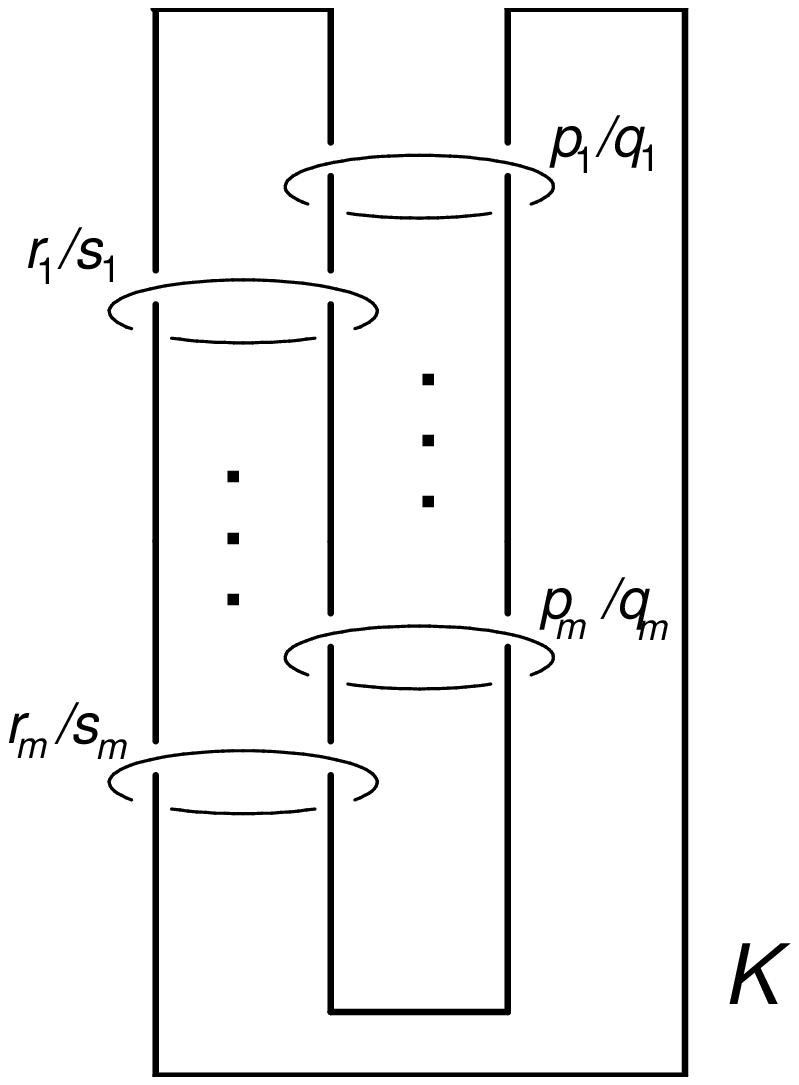}
\end{center}
\caption{} \label{Fig. 0}
\end{figure}

\begin{figure}[h]
\begin{center}
\includegraphics*[totalheight=16cm]{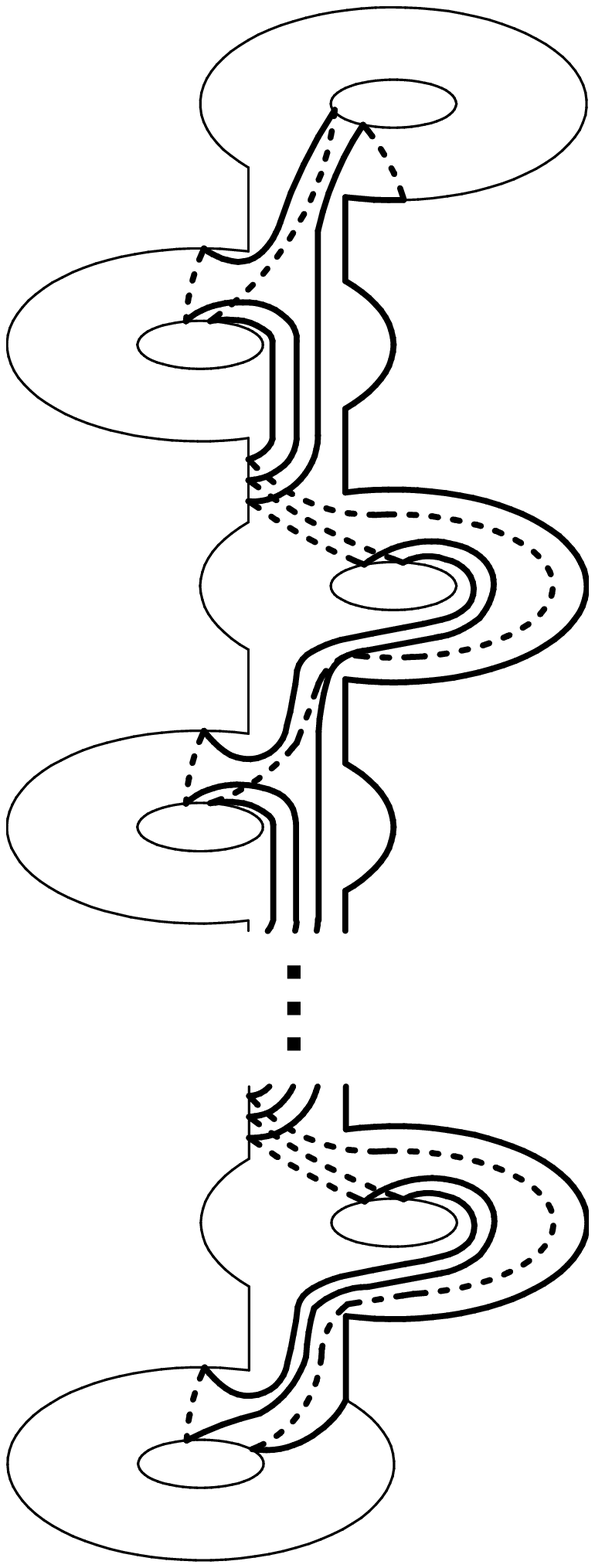}
\end{center}
\caption{} \label{Fig. 01}
\end{figure}

Take $2m-1$ ``vertical'' arcs $\sigma_1,\cdots,\sigma_{2m-1}$
consecutively connecting the components of $L$ and consider a
regular neighborhood $\H$ of
\hbox{$L\cup\sigma_1\cup\cdots\cup\sigma_{2m-1}$.} It is easy to
see that $\partial \H$ is a genus $2m$ Heegaard surface of $N$,
and $K$ can be embedded on it, as depicted in Figure \ref{Fig.
01}. Therefore $K$ is $(2m,0)$-knot and hence, as we already
noted, a $(2m,1)$-knot.
\end{proof}

Now we explicitly give Heegaard diagrams of periodic Takahashi
manifolds $T_n(p/q,r/s)$, for $p/q,r/s\ge 0$ (possibly equal to
$+\infty=1/0$). The cases $p/q< 0$ and/or $r/s< 0$ can be obtained
in a similar way.

\begin{figure}
\begin{center}
\includegraphics*[totalheight=18cm]{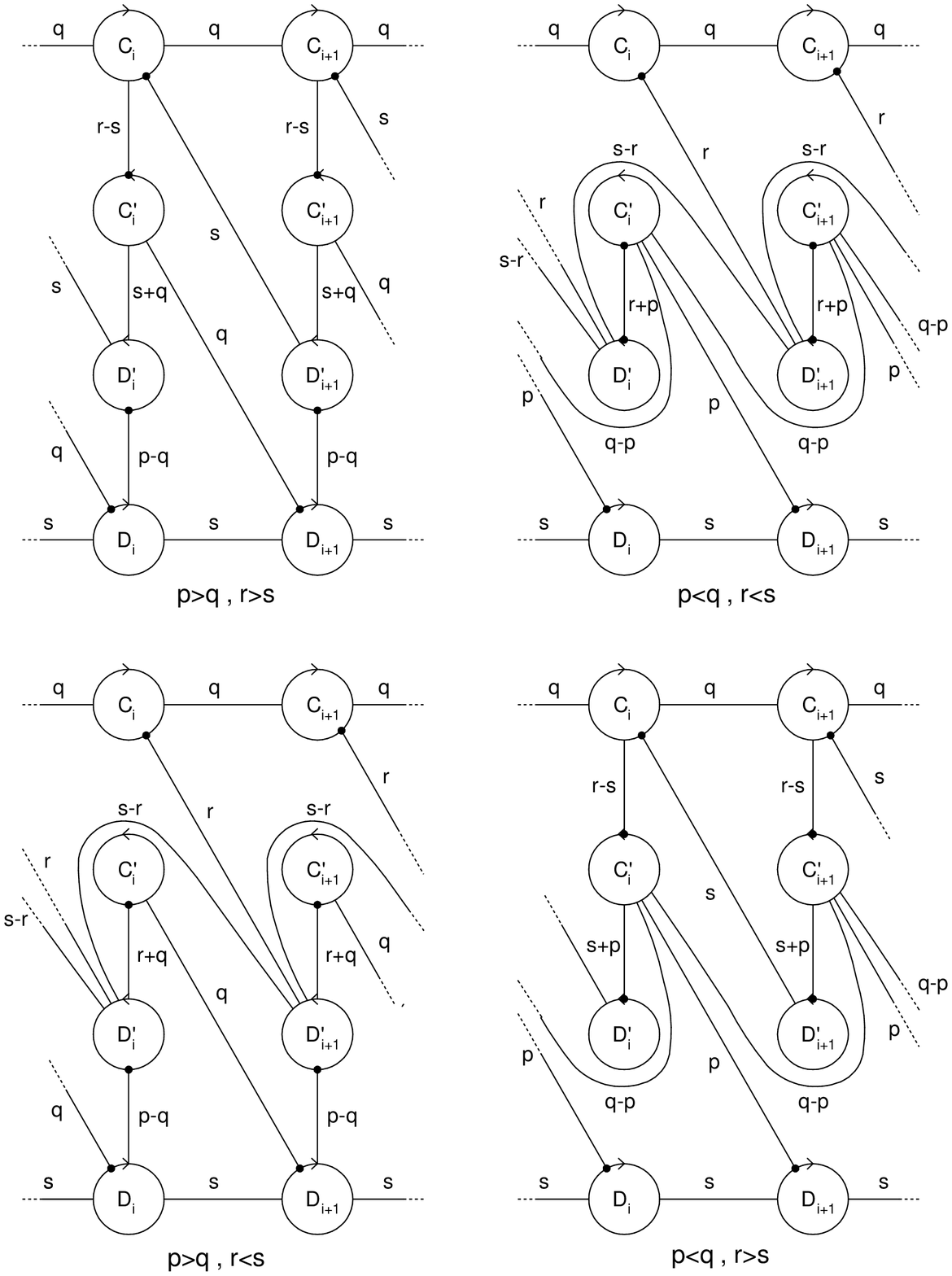}
\end{center}
\caption{Heegaard diagrams for $T_n(p/q,r/s)$} \label{Fig. 3}
\end{figure}

In Figure \ref{Fig. 3} an arc labelled $k$ denotes $k$ parallel
arcs. Furthermore, each diagram in the figure admits a cyclic
symmetry of order $n$; for $i=1,\ldots,n$, we identify $C_i$
(resp. $D_i$) with $C'_i$ (resp. $D'_i$) by gluing together the
first vertices of the marked ones, according to the orientations.

\begin{proposition}
The periodic Takahashi manifold $T_n(p/q,r/s)$, for $p/q,r/s\ge 0$
(possibly equal to $+\infty=1/0$), has Heegaard diagram as
presented in Figure \ref{Fig. 3}.
\end{proposition}

\begin{figure}
\begin{center}
\includegraphics*[totalheight=5.4cm]{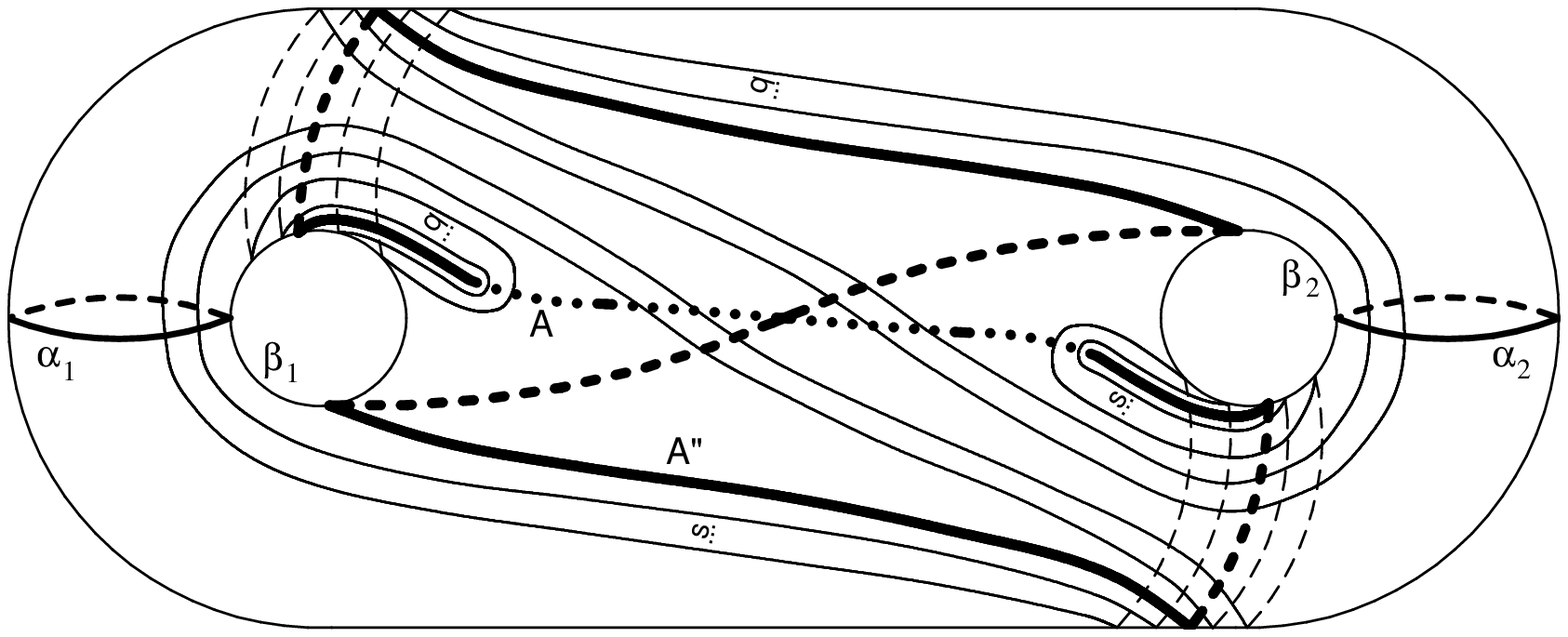}
\end{center}
\caption{} \label{Fig. 4}
\end{figure}

\begin{proof} Without loss of generality, we can assume that $p,q,r,s
\ge 0$. Let $\H_2$ be a genus two handlebody canonically embedded
in $\R^3$, as in Figure~\ref{Fig. 4}, which arises from Figure
\ref{Fig. 01} for $m=1$. Let us consider the Heegaard diagram of
genus two of $N=L(p,q)\# L(r,s)$ obtained by cutting the genus two
surface $\partial \H_2$ along the curves $\a_1$ and $\a_2$ and
gluing them back after a twist of $r/s$ and $p/q$ turns,
respectively, when $s,q\ne 0$. After this operation, the $s$
(resp. $q$) ``parallel'' thin curves in $\partial \H_2$ become the
curve $\f (\beta_1')$ (resp. $\f (\beta_2')$) of Theorem
\ref{basic}. The arcs $A''\subset\partial \H_2$ and $A\subset
\H_2$ (depicted in the figure by dashed lines) defines a knot
$A\cup A''$ which is exactly the knot $K$ of Proposition
\ref{Takahashi} for $m=1$. If $s=0$ (resp. $q=0$) the curve $\f
(\beta_1')$ (resp. $\f (\beta_2')$) is $\a_1$ (resp. $\a_2$).

By cutting $\H_2$ along $\beta_1$ and $\beta_2$ we obtain an open
Heegaard diagram of genus two of $N$. From Corollary
\ref{homology} we have $H_1(N-K) = \langle \a_1, \a_2, \bb \, | \,
p \a_2, r \a_1 \rangle$. So we have $z_1=\o_f(\a_1)=0$ and
$z_2=\o_f(\a_2)=0$ for each $n>1$. By performing the $n$-fold
strongly-cyclic covering determined by this monodromy, we obtain
the Heegaard diagrams of Figure \ref{Fig. 3}.
\end{proof}

\section{Fundamental groups of strongly-cyclic branched coverings} 

It is shown in \cite{M} that strongly-cyclic branched coverings of
$(1,1)$-knots admit geometric (i.e. induced by Heegaard diagrams)
cyclic presentations. In order to obtain a similar result for any
knot, we give the following definition.

\smallskip

\noindent {\bf Definition.} A balanced $mn$-generator presentation
of a group $\langle X\mid R\,\rangle$ with generator set $X=
\{x_{i,j}\mid 1\le i\le m\,,\,  1\le j\le n\}$ and relator set
\hbox{$R= \{r_{i,j}\mid 1\le i\le m\,,\,  1\le j\le n\}$} is said
to be an {\it $m$-word cyclic presentation\/} if there exist $m$
words $w_1, \ldots, w_m$ in the free group $F_{mn}$ generated by
$X$, such that the defining relations are of the form $r_{i,j} =
\theta_n^{j-1} (w_i)$, with $j=1,\ldots,n$, for each $i=1, \ldots,
m$, where \hbox{$\theta_n :F_{mn} \to F_{mn}$} is the automorphism
defined by $\theta_n (x_{i, j}) = x_{i, j+1}$ (with second
subscript mod $n$), for each $j=1,\ldots,n$ and $i=1, \ldots, m$.

\smallskip

We denote this group by $G_n(w_1, \ldots, w_m)$. For $m=1$ we have
the classical case of cyclically presented group (see \cite{Jo}).

There is a strict connection between strongly-cyclic branched
coverings of $(g,1)$-knots and $g$-word cyclic presentations, as
stated by the following:


\begin{theorem}\label{cyclic symmetry}
Every  $n$-fold strongly-cyclic branched covering $M$ of a
$(g,1)$-knot $K\subset N$ admits a Heegard splitting of genus $gn$
with cyclic symmetry of order $n$, which induces a $g$-word cyclic
presentation for the fundamental group of $M$.
\end{theorem}

\begin{proof} Since $K$ is a $(g,1)$-knot in $N$, we have $(N,K)=(\H_g,A)\cup_{\f}(\H_g',A')$, where
$\f:(\partial \H_g', \partial A')\to(\partial \H_g,\partial A)$ is
an attaching homeomorphism. Let $f:(M,f^{-1}(K))\to (N,K)$ be an
$n$-fold strongly-cyclic branched covering of $K$. Then
$\H_{gn}=f^{-1}(\H_g)$ and $\H_{gn}'=f^{-1}(\H_g')$ are both
handlebodies of genus $gn$. Moreover, $f^{-1}(A)$ and $f^{-1}(A')$
are both properly embedded arcs in $\H_{gn}$ and $\H_{gn}'$
respectively. Thus, we get a genus $gn$ Heegaard splitting
$(M,f^{-1}(K))=(\H_{gn},f^{-1}(A))\cup_{\F}(\H_{gn}',f^{-1}(A'))$,
where $\F:\partial \H_{gn}'\to\partial \H_{gn}$ is the lift of
$\f$ with respect to $f$. The group $\pi_1(\H_g,*)$ is generated
by the set $\{\la_1,\ldots,\la_g\}$ (see the notations introduced
before Proposition \ref{fundamental}). Let $\o_f$ be the monodromy
associated to the covering $f$, we can suppose that
$\o_f(\gamma)=1$. Moreover, we can assume that $\o_f (\la_i)=0$
for each $i=1,\ldots,g$. In fact, if $\o_{f}(\la_i)=k_i$ we
replace $\la_i$ with a loop homotopic to $\la_i\gamma ^{-k_i}$.
Taking a point $Q\in A$, let $\d$ be an arc from $*$ to $Q$ such
that $\d\cap A=Q$. Then $f^{-1}(Q)$ is a single point $\tilde Q\in
f^{-1}(A)$ and $f^{-1}(*)$ consists of $n$ points $\tilde
*_1,\ldots,\tilde *_n$. For $j=1,\ldots,n$ and $i=1,\ldots,g$, let
$\tilde\aa_{ij}$ and $\tilde\d_j$ be the lifts (with respect to
$f$) of $\la_i$ and $\d$ respectively, both containing $\tilde
*_j$. Then, for each $i=1,\ldots,g$, the $n$ loops
$$\aa_{i1}=\tilde\d_1^{-1}\cdot\tilde\aa_{i1}\cdot\tilde\d_1,
\, \ldots, \,\, \aa_{in}=\tilde\d_n^{-1} \cdot \tilde\aa_{in}
\cdot \tilde\d_n$$ are cyclically permutated by a generator $\Psi$
of the group of covering transformations and the set
$\{\aa_{ij}\mid 1\le i\le g\,, \,  1\le j\le n\}$ generates
$\pi_1(\H_{gn},\tilde Q)$. Let $\{E_1',\ldots,E_g'\}$ be a
complete system of meridian disks for $\H_g'$ such that $E'_i\cap
A'=\emptyset$, for each $i=1,\ldots,g$. Then  $f^{-1}(E_i')$
consists of $n$ disks $\tilde E_{i1}',\ldots,\tilde E_{in}'$,
which are cyclically permutated by $\Psi$, for each
$i=1,\ldots,g$, and $f^{-1}(E_1'\cup\cdots\cup E_g')$ is a
complete system of meridian disks for the handlebody $\H_{gn}'$.
The curves $\F(\partial\tilde E_{11}'),\ldots,\F(\partial\tilde
E_{gn}')$ give the relators for the presentation of
$\pi_1(M,\tilde Q)$ induced by the Heegaard splitting. Since, for
each $i=1,\ldots,g$, the generators $\aa_{i1},\ldots,\aa_{in}$ and
the relator curves $\partial\tilde E_{i1}',\ldots,\partial\tilde
E_{in}'$ are cyclically permutated by $\Psi$, each $\partial\tilde
E_{ij}'$ corresponds to a cyclic permutation of the generators in
the word $w_i$.
\end{proof}

\medskip

Now we describe an algorithm to obtain a presentation of the
fundamental group of a strongly-cyclic branched covering of a
$(g,1)$-knot.

Let $\o_{f}$ be the monodromy of an $n$-fold strongly-cyclic
branched covering of a $(g,1)$-knot $K$ in a 3-manifold $N$ and
let us consider the presentation of the knot group obtained in
Proposition \ref{fundamental}. Following the proof of Theorem
\ref{cyclic symmetry}, we choose, for each $i=1,\ldots,g$, a new
generator $\wa_i =\la_i\gamma^{-\o_{f}(\bar\a_i)}$ and we get
$\p_1(N-K,*) = \langle \wa_1,\ldots,\wa_g, \gamma\,|\,\bar
r_i(\wa_1, \ldots, \wa_g,\gamma), \quad i=1, \ldots, g \rangle$,
with $\bar r_i(\wa_1, \ldots, \wa_g,\gamma) =
r_i(\wa_1\gamma^{\o_{f}(\bar\a_1)}, \ldots, \wa_g
\gamma^{\o_{f}(\bar\a_g)}, \gamma)$.

Suppose now that, for each $i=1,\ldots,g$, we have
$$r_i(\wa_1,\ldots,\wa_g,\gamma)=\wa_{j_1}^{\ee_{ij_1}}\gamma^{\eta_{ij_1}}
\wa_{j_2}^{\ee_{ij_2}}\gamma^{\eta_{ij_2}}\ldots\wa_{j_{\l_i}}^{\ee_{ij_{\l_i}}}\gamma^{\eta_{ij_{\l_i}}}$$
where $\ee_{ij_l}\ne 0$ for each $l=1,\ldots,\l_i$ and
$j_1,\ldots,j_{\l_i}$ are not necessary distinct.

Then we have the following:

\begin{proposition} \label{words} The fundamental group of
the $n$-fold strongly-cyclic branched covering of a $(g,1)$-knot
$K$, with monodromy $\o_f$, admits the $g$-word cyclic
presentation $G_n(w_1,\ldots,w_g)$, where, for each
$i=1,\ldots,g$,
$$w_i=x_{{j_1},1}^{\ee_{ij_1}}x_{{j_2},k_{i2}}^{\ee_{ij_2}}\cdots x_{j_{\l_i},k_{i\l_i}}^{\ee_{ij_{\l_i}}},$$
with $k_{il}\equiv 1+\sum_{h=1}^{l-1}\eta_{ij_h} \mod n$, for
$l=2,\ldots,\l_i$.
\end{proposition}
\begin{proof}
Let $M$ be the $n$-fold strongly-cyclic branched covering of $K$
determined by $\o_f$. From the proof of Theorem \ref{cyclic
symmetry},  the fundamental group of $M$ admits the set of
generators $\{\a_{ij}\mid 1\le i\le g\,, \,  1\le j\le n\}$, where
$\a_{i1},\ldots,\a_{in}$ are the components  of the lifting of
$\wa_i$. The relators are (homotopic to) $\F(\partial \tilde
E'_{11}),\ldots, \F(\partial \tilde E'_{gn})$, where $\F$ is the
lifting of $\ps$ with respect to $f$, and each $\partial \tilde
E'_{ij}$ is a component of the lifting of a meridian disc $E'_i$
of $\H_g'$ such that $E'_i\cap A'=\emptyset$. For each
$i=1,\ldots,g$, we can choose $E'_i$ such that $\partial
E'_i=\tau(\b_i)$, and therefore the relators are (homotopic to)
the components of the lifting of $r_i(\la_1,\ldots,\la_g,\gamma)$,
or equivalently $\bar r_i(\wa_1,\ldots,\wa_g,\gamma)$
($i=1,\ldots,g$), which arise from the relations
$\lb'_i=\ps(\lb'_i)$. Now, since $\o_{f}(\bb)=1$, a factor
$\gamma^k$ lifts to a path connecting the point of $f^{-1}(*)$ in
the sheet $j$ with the corresponding point in the sheet $j+k$ (mod
$n$). The result immediately follows.
\end{proof}


\noindent {\bf Example.} Let us consider the $(2,1)$-knot $K=A\cup
A''$ of Figure \ref{Fig. 4}, with $q=s=1$. In this case $K\subset
N=L(p,1)\#L(r,1)$, and $$\pi_1(N-K)=\langle \a_1,\a_2,\g\mid
\a_2\a_1^{-r}\g\a_2^{-1}\g^{-1},\a_1\g\a_2^p\a_1^{-1}\g^{-1}
\rangle .$$ By applying Proposition \ref{words}, with
$\o_f(\a_1)=\o_f(\a_2)=0$, we get
$\pi_1(T_n(p/1,r/1))=G_n(w_1,w_2)$, where
$w_1=x_{2,1}x_{1,1}^{-r}x_{2,2}^{-1}$ and
\hbox{$w_2=x_{1,1}x_{2,2}^px_{1,2}^{-1}$,} which is exactly the
group presentation obtained in \cite{Ta}.

\end{document}